\documentclass[12pt]{article}
\usepackage{fullpage}
\usepackage{verbatim}
\usepackage{graphics, graphicx}
\usepackage{amsthm, amsmath, amssymb}

\renewcommand{\comment}[1]{ }
\newcommand{\nn}{\nonumber\\}
\newcommand{\ds}{\displaystyle}

\def\NW{{\sf NW}}\def\NN{{\sf N}}
\def\SE{{\sf SE}}
\def\NE{{\sf NE}}

\def\seta{\mathcal{S}}
\def\setb{\mathcal{T}}

\def\y{Y}
\def\x{X}
\def\oy{\overline{\y}}

\def\bY{A}
\def\bX{B}

\def\obY{\overline{\bY}}
\def\obX{\overline{\bX}}

\def\xn{\underline{x}}
\newcommand\bZ{\ensuremath{\mathbb Z}}
\newcommand\bC{\ensuremath{\mathbb C}}

\newtheorem{thm}{Theorem}[section]
\newtheorem{cor}[thm]{Corollary}

\newtheorem{lemma}[thm]{Lemma}

\numberwithin{equation}{section}

\newtheoremstyle{marnidef}{\topsep}{\topsep}%
     {}
     {}
     {\bfseries}
     {.  }
     {.5em}
     {\thmname{#1}\thmnumber{#2 }\thmnote{{\em #3}}}
\theoremstyle{marnidef}
\newtheorem{defn}{Definition }[section]

\begin{document}
\title{Two Non-holonomic Lattice Walks in the Quarter Plane}%
\author{
\begin{minipage}{8cm}\begin{center}
   Marni Mishna\\ 
   {\small Dept. of Mathematics}\\ 
   {\small Simon Fraser University}\\
   {\small Burnaby, Canada}
\end{center}
\end{minipage}
\begin{minipage}{8cm}
\begin{center}
   Andrew Rechnitzer\\
   {\small Dept. of Mathematics}\\ 
   {\small University of British Columbia}\\
   {\small Vancouver, Canada}
\end{center}
\end{minipage}
}
\date{\mbox{}}
\maketitle
\begin{abstract}
We present two classes of random walks restricted to
the quarter plane whose generating function is not holonomic. The
non-holonomy is established using the iterated kernel method, a recent
variant of the kernel method. This adds evidence to a recent
conjecture on combinatorial properties of walks with holonomic
generating functions. The method also yields an asymptotic expression
for the number of walks of length $n$. 

\small {\bfseries keywords:} 
random walks, enumeration, holonomic, generating functions, kernel
method
\end{abstract}

\section{Introduction}
Previous studies of random walks on regular lattices had led some to
conjecture that random walks in the quarter plane should all have
holonomic generating functions. This is not unreasonable, given that
walks confined to the half plane have algebraic generating functions
(\cite{BaFl02} proves this for directed paths, but their results
extend to random walks). This was, however, disproved by
Bousquet-M\'elou and Petkov\v sec with their proof that knight's walks
confined to the quarter plane are not holonomic~\cite{BoPe03}.

In this paper we give two new examples of random walks whose
generating functions are not holonomic. We use a technique developed
in a recent study of self-avoiding walks in wedges~\cite{JaPeRe06}
called the iterated kernel method. This is an adaptation of the kernel
method~\cite{MBM05,Prodinger03}, in which we express a generating
function in terms of iterates of a kernel solution.  We can show that
this generating function has an infinite number of poles, which is
incompatible with holonomy.

The aim of this work is two-fold. The first goal is to illustrate a
potentially general technique to prove non-holonomy. Secondly, this
work was completed in the context of a general generating function
classification of all nearest neighbour lattice walks~\cite{MishnaXX},
and lends evidence to a general conjecture on the connection between
different symmetries of the step sets and the analytic nature of their
generating functions.

\subsection{Walks and their generating functions} 
\begin{figure}[t]
\label{fig:sample}
\center
\includegraphics[width=6cm]{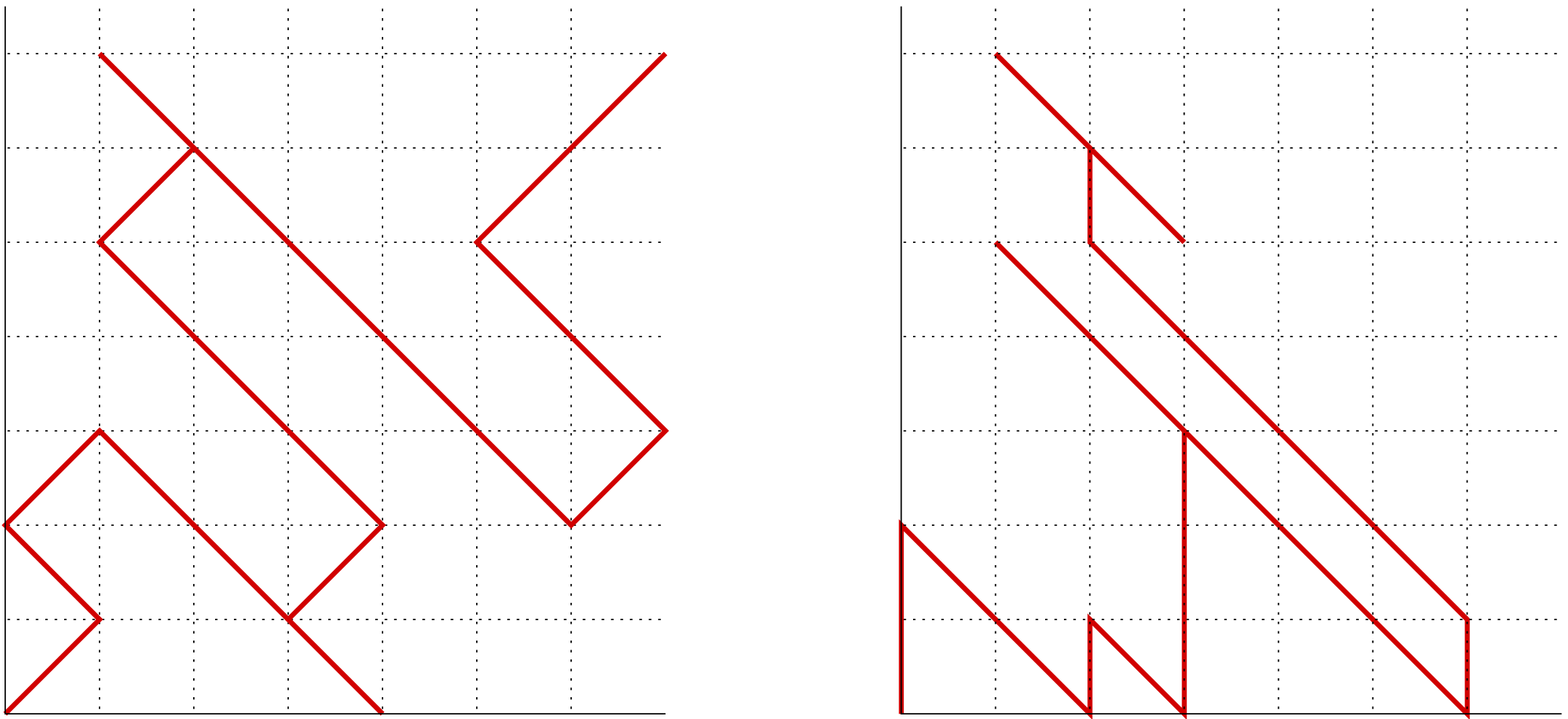}
\caption{Sample walks with steps from $\seta=\{\NW, \NE, \SE\}$ and $\setb=\{\NW,\NN, \SE\}$} 
\end{figure}

The objects under consideration are walks on the integer lattice
restricted to the quarter plane with steps taken from $\seta=\{(-1,1),
(1,1), (1,-1) \}$ in the first case, and $\setb=\{(-1,1), (0,1),
(1,-1) \}$ in the second case. We use the compass notation and label the directions of these two sets as $\seta=\{\NW, \NE, \SE\}$ and
$\setb=\{\NW, \NN, \SE\}$. Two sample walks are given in
Figure~\ref{fig:sample}. 

We associate to these two steps sets two power series each:~$W(t)$ a
counting (ordinary, univariate) generating function and~$Q(x,y;t)$ a
complete (multivariate) generating function. The series~$W(t)$, is
the ordinary generating function for the number of walks. It is a
formal power series where the coefficient of $t^n$ is the number of
walks of length~$n$. The complete generating function $Q(x,y;t)$
encodes more information. The coefficient of $x^iy^jt^n$
in~$Q(x,y;t)$ is the number of walks of length $n$ ending at the
point~$(i,j)$. Remark that the specialization $x=y=1$ in the complete
generating function is precisely the counting series,
i.e. $Q(1,1;t)=W(t)$. If the choice of step set is not clear, we add
an index. 

In part, our interest in the complete generating function stems from
the fact it satisfies a very useful functional equation which we derive using the
recursive definition of a walk: a walk of length~$n$ is a
walk of length $n-1$ plus a step. The quarter plane condition asserts
itself by restricting our choice of step should the smaller walk end
on an axis.

The step set $\seta$ leads to the following equation:
\begin{equation}\label{eq:seta}
  Q(x,y;t) = 1 + t \left(xy+\frac{x}{y}+\frac{y}{x}\right) Q(x,y;t) 
  -t \frac{x}{y}Q(x,0;t)  -t \frac{y}{x}Q(0,y;t),
\end{equation}
and set $\setb$ defines a similar equation:
\begin{equation}\label{eq:setb}
  Q(x,y;t) = 1 + t \left(y+\frac{x}{y}+\frac{y}{x}\right) Q(x,y;t) 
  -t \frac{x}{y}Q(x,0;t)  -t \frac{y}{x}Q(0,y;t).
\end{equation}

\subsection{Properties of holonomic functions} 
We are interested in understanding the analytic nature of the
generating functions. This nature gives a basic first classification
of structures, and also some general properties such as general form
of the asymptotic growth of the coefficients. See for example
Bousquet-M\'elou's recent work on classifying combinatorial families with
rational and algebraic generating functions~\cite{MBM06}. We are interested in generating functions which are
holonomic, also known as D-finite. Let $\xn=x_1, x_2,
\dots, x_n$.

\begin{defn}[holonomic function]
A multivariate function~$G(\xn)$ is {\em holonomic\/} if the
vector space generated by the partial derivatives of~$G$ (and their
iterates), over rational functions of $\xn$ is finite
dimensional. This is equivalent to the existence of $n$ partial
differential equations of the form%
\[p_{0,i} f(\xn)+ p_{1,i} \frac{\partial f(\xn)}{\partial x_i}+\ldots 
   + p_{d_i,i} \frac{\partial^{d_i}f(\xn)}{(\partial x_i)^{d_i}}=0,
\] for $i$ satisfying $1\leq i\leq n$,
and where the $p_{j,i}$ are all polynomials in $\xn$. 
\end{defn}
The basic feature of these functions that we use is the fact that they
have a finite number of singularities. For example, in the univariate
case, these are poles which appear are zeroes of~$p_{d_1,1}$. 
\subsection{Main results}
We show that both the complete and counting
generating functions of the walks with steps from $\seta=\{\NW,\NE,\SE\}$ and
$\setb=\{\NW,\NN,\SE\}$ are not holonomic using this singularity
property of holonomic functions. In fact, it suffices
to show that their counting generating function $W(t)$ is not
holonomic, as holonomic functions are closed under algebraic
substitution~\cite{Lipshitz89,Stanley80}, and we have $W(t)=Q(1,1;t)$. We show that $W(t)$
has an infinite number of poles and thus is not holonomic~\cite{Stanley80}.
This will give the first main theorem that we present.

\begin{thm}\label{thm:notDfin}
Neither the complete generating function~$Q_{\seta}(x,y;t)$, nor the counting
generating function $W_{\seta}(t)=Q_{\seta}(1,1;t)$ of nearest-neighbour
walks in the first quadrant with steps from $\{\NE,\SE,\NW\}$ 
are holonomic functions with respect to their variable sets.
\end{thm}

To prove this, we use the {\em iterated kernel method}, as described
by Janse van Rensburg {\em et al.\/} in~\cite{JaPeRe06} in counting
partially directed walks confined to a wedge.

The basic idea is similar to other variants of the kernel method in
that we write our fundamental equation in the kernel format, determine
particular values of~$x$ and~$y$ which fix the kernel, and then
generate new relations. For these walks, the kernel iterates form an
infinite group whose elements each introduce new poles into the generating
function. The main difficulty in proving that these generating functions
are not holonomic is proving that the singularities do not cancel.

As step set~$\setb$ is not symmetric in the line $x=y$, the relations
obtained are more complex and while the argument is essentially the
same, the details of the proof of its non-holonomy is more
complicated, as we see in Section~\ref{sec:setb}. This will give the
second main theorem.

\begin{thm}\label{thm:Q11}
Neither the complete generating function~$Q_{\setb}(x,y;t)$, nor the counting
generating function $W_{\setb}(t)=Q_{\setb}(1,1;t)$ of nearest-neighbour
walks in the first quadrant with steps from $\{\NW,\NN,\SE\}$ 
are holonomic functions with respect to their variable sets.
\end{thm}
 
\section{The non-holonomy of  $Q_\seta(x,y;t)$}
\label{sec:seta}
We use the iterated kernel method to find a closed form expression for
the generating function $Q_\seta(x,y;t)$. We then analyse this
expression to prove this function, and the counting generating
function $W_\seta(t)$ are not holonomic.
\subsection{Defining the iterates of the kernel solutions}
To begin, consider the kernel version of Equation~\eqref{eq:seta}:
\begin{equation}\label{eqn:kform}
 \left(xy - tx^2y^2 - tx^2 - ty^2\right) Q(x,y) = xy - tx^2 Q(x,0) - ty^2Q(y,0).
\end{equation}
Here, for brevity we write $Q(x,y;t)$ as $Q(x,y)$, and have
used the $x \leftrightarrow y$ symmetry to rewrite $Q(0,y)$ as $Q(y,0)$.

Since the kernel~$K(x,y)=xy - tx^2y^2 - tx^2 - ty^2$ is a quadratic,
it has two solutions as a function of $y$
\begin{equation}
  \y_{\pm 1}(x) = \frac{x}{2t(1+x^2)} \left(1 \mp \sqrt{1-4t^2(1+x^2)} \right).
\end{equation}
Note that as formal power series in $t$, these roots are:
\begin{align*}
  \y_{+1}(x) & = xt + O(t^3) \\
  \y_{-1}(x) & = \frac{x}{(1+x^2)} \frac{1}{t} - xt + O(t^3).\\
\end{align*}
Further, $\y_{+1}(x)$ is a power series in $t$ with polynomial
coefficients in $x$. We will need to form compositions of these roots.
Write $\y_n(x) = \left( \y_1 \circ \right)^n(x)$, and likewise for
$\y_{-n}$.

\begin{lemma}\label{lem kern roots seta}
 The kernel roots obey the relation:
 \begin{equation}
    \label{eqn:y_inv}
    \y_{+1} ( \y_{-1} (x) ) = x\quad \text{ and }\quad \y_{-1} ( \y_{+1} (x) ) = x.\end{equation}
  and the set $\{\y_n|n\in\bZ\}$ forms a group, under the operation
  $\y_n(\y_m(x))=\y_n\circ \y_m=\y_{n+m}$, with identity~$\y_0=x$.
Furthermore, they also obey the relation:
  \begin{equation}
    \label{eqn:first_recip}
    \frac{1}{\y_1(x)} + \frac{1}{\y_{-1}(x)} = \frac{1}{tx},
  \end{equation}
  which extends (by substituting $x = \y_{n-1}(x)$) to:
 \begin{equation}\label{eqn:recip}
   \frac{1}{\y_n(x)} = \frac{1}{t \y_{n-1} } - \frac{1}{\y_{n-2}(x)}.
 \end{equation}
\end{lemma}
\begin{proof}

  To prove Equation~\eqref{eqn:y_inv}, consider the four compositions of the kernel roots
  $\y_{\pm}(\y_{\pm}(x))$ and denote them by $\{C_i\}$. We form
  the product \mbox{$\ds \prod_{i=1}^4 (C-C_i)$}, for formal parameter
  $C$, and one may verify:
  \begin{displaymath}
    (x-C)^2 \left( (t^2+x^2) C - x (1-2 t^2) C+t^2 x^2 \right) = 0.
  \end{displaymath}
  Hence two of the compositions are the identity.
  
  There are $\binom{4}{2}$ possibilities, and we can exclude 4 of
  these by consistency arguments --- \emph{i.e.} they will imply that
  $\y_{+} = \y_{-}$ which is not true. This leaves us to check that
  either $\y_+ \circ \y_- = \y_-\circ \y_+ = x$ or $\y_+ \circ \y_+ = \y_-\circ
  \y_- = x$. A manual calculation then shows that $\y_{-} \circ
  \y_{+} = \y_{+} \circ \y_{-} = x$.

  Equation~\eqref{eqn:first_recip} follows since $\y_{\pm 1}$ are
  roots of a quadratic --- in particular $\y_{+1} \y_{-1}$ and
  $\y_{+1} + \y_{-1}$ are rational functions of $x$ and $t$.
\end{proof}
\comment{
     \begin{itemize}
     \item If $\y_+ \circ \y_+ = \y_+ \circ \y_- = x$ then $\y_+ \circ
       \y_+ \circ \y_+ = \y_+ = \y_+ \circ \y_+ \circ \y_- = \y_-$
       which is a contradiction.
     \item Similarly $\y_+ \circ \y_+ = \y_- \circ \y_+ = x$, $\y_-
       \circ \y_- = \y_+ \circ \y_- = x$ and $\y_- \circ \y_- = \y_-
       \circ \y_+ = x$ also lead to contradictions.
     \item If $\y_+ \circ \y_- = \y_- \y_+ = x$ or if $\y_+ \circ \y_+
       = \y_- \y_- = x$ then we cannot argue on consistency grounds.
     \end{itemize}

     We now show that $\y_{-1} \circ \y_{+1} = x$. Writing $\Delta =
     1-4t^2(1+x^2)$ we have:
     \begin{align}
       \y_{-1} (\y_{+1}(x)) & = x \frac{\left(1+2x^2-\sqrt{\Delta}
         \right) \left(1+ x^2 + \sqrt{x^4 + 2x^2\sqrt{\Delta} + \Delta}
         \right) }{4 (1+x^2)(t^2+x^2)} \nn & = x
       \frac{\left(1+2x^2-\sqrt{\Delta} \right) \left(1+ x^2 + (x^2 +
           \sqrt{\Delta}) \right) }{4 (1+x^2)(t^2+x^2)} \qquad \mbox{is
         there a sign change?} \nn & = x \frac{(1+2x^2)^2-\Delta}{4
         (1+x^2)(t^2+x^2)} \nn & = x
       \frac{x^2+t^2+x^4+t^4}{(1+x^2)(t^2+x^2)} = x
     \end{align}
     One may similarly show that $\y_{+} \circ \y_{-} = x$, and that
     the other two compositions are not equal to $x$.

}

\subsection{An expression for $Q(x,0)$ in terms of the iterates}
We are now able to find a formal expression for the generating function.
\begin{thm}
  \label{thm Q11 symm}
  The generating function $Q(x,0)$ is given by
  \begin{equation}\label{eqn:Qxo}
    Q(x,0) = \frac{1}{x^2 t} \sum_{n \geq 0} (-1)^n \y_n(x) \y_{n+1}(x).
  \end{equation}
  and so $Q(1,1)$ is given by
  \begin{equation}\label{eq:infsing}
    Q(1,1) = \frac{1-2tQ(1,0)}{1-3t}
    = \frac{1 - 2\sum_{n \geq 0} (-1)^n \y_n(1) \y_{n+1}(1)}{1-3t}.
\end{equation}
\end{thm}
\begin{proof}
  By construction we have $K(x,\y_1(x)) = 0$, thus substituting $y =
  \y_1(x)$ into Eq.~\eqref{eqn:kform} gives (after a little tidying)
  \begin{equation*}
    Q(x,0) = \frac{\y_1(x)}{x} \frac{1}{t} - \frac{\y_1(x)^2}{x^2} Q(\y_1(x), 0).
  \end{equation*}
  Now we substitute~$x = \y_{n}(x)$ into this equation to obtain
  \begin{equation*}
    Q(\y_n(x),0) = \left(\frac{\y_{n+1}(x)}{\y_n(x)}\right) \frac{1}{t} 
    - \left(\frac{\y_{n+1}(x)}{\y_{n}(x)}\right)^2 Q(\y_{n+1}(x), 0). 
  \end{equation*}
  Using this expression for~$Q(\y_n(x),0)$ for increasing~$n$, we can
  iteratively generate a new expression for~$Q(x,0)$:
  \begin{align*}
    Q(x,0) 
    &= \frac{\y_1(x)}{x} \frac{1}{t} 
    - \frac{\y_1(x)^2}{x^2} \left(
      \frac{\y_2(x)}{\y_1(x)} \frac{1}{t} - \frac{\y_2(x)^2}{\y_1(x)^2} Q(\y_2(x), 0)
    \right) \nn
    &= \frac{\y_0(x)\y_1(x)}{x^2t}- \frac{\y_1(x) \y_2(x)}{x^2t}
    + \left( 
      \frac{\y_2(x)}{x} \right)^2 Q(\y_2(x),0)\nn
    &=\dots\nn
    &= \frac{1}{x^2 t} \sum_{n=0}^{N-1} (-1)^n\y_n(x) \y_{n+1}(x)
    + (-1)^N \left(\frac{\y_N(x)}{x} \right)^2 Q(\y_N(x),0).
\end{align*}
Since $\y_1(x) = xt + O(t^3)$ is a power series in $t$ with positive polynomial
coefficients in $x$, it follows that $\y_n(x) = xt^n + o(xt^n)$. Hence
we have that $\ds \lim_{N \to \infty} \y_N(x) = 0$ as a formal power
series in $t$ and the theorem follows.
\end{proof}

\subsection{The abundant singularities of  $W(t)$}
\label{sec:sings}
Our path to proving Theorem~\ref{thm:notDfin} is now clear. We demonstrate  an
infinite set of singularities for $W(t)=Q(1,1)$ coming from the set of
non-conflicting poles of~$\y_n(1)$. Given equation~\eqref{eq:infsing},
there are 3 potential sources of singularities for $Q(1,1)$:
\begin{enumerate}
\item the simple pole at $t = 1/3$;
\item singularities from the $\y_n(1;t)$;
\item divergence of the sum.
\end{enumerate}
We now show that the singularity at $t=1/3$ is the dominant
singularity and that there are an infinite number of other
singularities (given by (2)), and that the series is convergent elsewhere.

\begin{lemma}
  \label{lem simple pole}
  The generating function $W_{\seta}(t)$ has a simple pole at
  $t=1/3$.
\end{lemma}
\begin{proof}
  To show that $t=1/3$ is a simple pole, we show that the numerator of
  equation~\eqref{eq:infsing}, $1-2tQ(1,0)$, is non-zero and
  absolutely convergent at this point.

  We show in the proof of Corollary~\ref{cor asympt1} that the radius
  of convergence of $Q(1,0)$ is bounded below by $8^{-1/2}\approx
  0.335$, and thus at $|t|=1/3$ the numerator converges absolutely.

  Next we show that the singularity is not removable by proving that
  $1-2tQ(1,0)$ is non-zero at $t=1/3$. Note that $\y_0(1;1/3) = 1$
  and $\y_1(1;1/3) = 1/2$. Using the recurrence~\eqref{eqn:recip}
  for~$\y_n$ we have
\begin{equation*}
  \frac{1}{\y_n(1;1/3)} = \frac{3}{ \y_{n-1}(1;1/3)} - \frac{1}{\y_{n-2}(1;1/3)}
\end{equation*}
which leads to $\y_n(1;1/3) = 1/F_{2n}$, \emph{i.e.} the reciprocal of
the~$2n^{\mathrm{th}}$ Fibonacci number. Hence $\y_n \y_{n+1} \to 0$
and so the alternating sum~$Q(1,1;1/3)$ is convergent. Further, we may
write
\begin{align}
  \sum_{n \geq 0} \frac{(-1)^n}{F_{2n} F_{2n+2}}
  & = \frac{1}{2} - \sum_{k \geq 1} \frac{F_{4k+2}-F_{4k-2}}{F_{4k-2}F_{4k} F_{4k+2}} \nn
  & = \frac{1}{2} - \frac{1}{10} + \sum_{k \geq 1} \frac{F_{4k+4}-F_{4k}}{F_{4k}F_{4k+2} F_{4k+4}}.
\end{align}
Hence the sum is bounded strictly between $2/5$ and $1/2$ and
$1-2Q(1,0;1/3)/3$ is non-zero, and thus $t=1/3$ is a simple pole of $Q(1,1)$.
\end{proof}

\subsection{The singularities are distinct}
Next, we show that this function possesses an infinite number of
singularities --- coming from simple poles of the $\y_n$. In order to
do this we make the substitution $t \mapsto \frac{q}{1+q^2}$ which
allows us to write $\y_n$ in closed form.

\begin{lemma}
  Define $\oy_n(q) := \left(
    \y_n\left(1;\frac{q}{1+q^2}\right) \right)^{-1}$. Then
  \begin{align}\label{eqn ynq}
    \overline{\y_n}(q) 
    & = \left(\frac{q-\overline{\y_1}(q) }{1-q^2}  \right) \cdot q^{1-n} +
    \left(\frac{1-q\overline{\y_1}(q)}{1-q^2} \right) \cdot q^n \nn
    & =   \left(\frac{1-q^{2n}}{1-q^2} \right) \cdot q^{1-n} \cdot  \overline\y_1(q)
    -  \left(\frac{1-q^{2n-2}}{1-q^2} \right) \cdot q^{2-n}.
  \end{align}
\end{lemma}
\begin{proof}
  Substitute $t = \frac{q}{1+q^2}$ into
  equation~\eqref{eqn:recip}. The roots of the characteristic equation
  are simply $q$ and $1/q$ and the recurrence can be solved using
  standard methods.
\end{proof}

Since we now have $\y_n$ in closed form we can determine some facts
about the location of its poles.
\begin{lemma}\label{thm:work}
Suppose $q_c$ is a zero of~$\overline\y_n(q):=\y_n\left(1;\frac{q}{1+q^2}\right)^{-1}$, and that $q_c\neq 0$. Then 
\begin{enumerate}
\item $q_c$ is a solution of $\left( q^{2n} + q^{-2n} \right) + \left(q^2 + q^{-2} \right) = 4$;
\item $q_c$ is not a $k^{\text{th}}$ root of unity for any $k$; and
\item for all $k \neq n$,  $\overline\y_k(q_c) \neq 0$.
\end{enumerate}
\end{lemma}

\begin{proof}
If $\overline\y_n(q_c) = 0$ we must have
\begin{equation*}
  q_c (1-q_c^{2n}) \overline\y_1(q_c)  = (q_c^2-q_c^{2n}).
\end{equation*}
Using the equation $\overline\y_1(q) =\left(1+q^2+\sqrt{1-6q^2+q^4} \right)/2q$, this can be rearranged to
give
\begin{subequations}
\begin{align*}
  (1-q_c^{2n}) \left( 1+q_c^2+\sqrt{1-6q_c^2+q_c^4} \right) & = 2 (q_c^2-q_c^{2n}),\\
  (1-q_c^{2n}) \sqrt{1-6q_c^2+q_c^4} & = 2 (q_c^2-q_c^{2n}) - (1-q_c^{2n})(1+q_c^2)\\
  & = - (1-q_c^2)(1+q_c^{2n}),\\
  \text{and so, }\qquad(1-q_c^{2n})^2 (1-6q_c^2+q_c^4) & = (1-q_c^2)^2 (1+q_c^{2n})^2.
\end{align*}
\end{subequations}
We expand and collect terms to arrive at the equation
\begin{equation*}
  q_c^{4n} + q_c^{2n} \left(\frac{1-4q_c^2+q_c^4}{q_c^2} \right) +1 = 0,
\end{equation*}
which can in turn be reduced to
\begin{equation*}
  \left( q_c^{2n} + q_c^{-2n} \right) + \left(q_c^2 + q_c^{-2} \right) = 4.
\end{equation*}

It is clear that $q = \pm 1$ are solutions of the above
equation. However, substituting $q = \pm 1$ into
equation~\eqref{eqn ynq} then we find $\y_1(\pm1) = \pm (1-i)/2$ and
\begin{align}
  \overline \y_n(1) &= n \overline \y_1(1) - (n-1)\nn
  \overline \y_n(-1) &= (-1)^{n+1} \left( n  \overline \y_1(-1) -(n-1) \right).
\end{align}
Hence neither equation can be satisfied and $q = \pm 1$ are not a
poles of $\y_n$. Now if $q_c = e^{i \theta}$, it follows that $\cos(2n
\theta) + \cos(2 \theta) = 2$ so that the only possible zeros on the
unit circle are at $q_c = \pm 1$, which we have already excluded. All
the zeros lie off the unit circle.

Next, we prove part~(3) of the lemma. First consider the possibility
that $\overline\y_{n+1}(q_c) = 0$. If this is true, since the sequence
of $\overline\y_n$ satisfies a homogeneous three-term recurrence, it
follows that $\overline\y_{n+k}(q_c) = 0$ for all $k \geq 0$. It also
follows that all $\overline\y_{n-k}(q_c)$ must be zero. However, we
can compute the zeros of $\overline\y_1(q)$ explicitly and the only
possible zero is at $q=0$. Hence either $q_c = 0$ or $\overline\y_{n
  \pm 1}(q_c) \neq 0$.

Using the recurrence $\overline\y_k(q_c) = \left(\frac{1+q_c^2}{q_c}\right)
\overline\y_{k-1}(q_c) - \overline\y_{k-2}(q_c)$, and the hypothesis that $\overline\y_n(q_c) = 0$,
we solve the recurrence for  the $\overline\y_{n+k}(q_c)$. This gives
\begin{equation}\label{eqn:bnq}
  \overline\y_{n+k}(q_c) = \overline\y_{n+1}(q_c) \frac{1-q_c^{2k}}{(1-q_c^2) q_c^{k-1}}.
\end{equation}
Since $\overline\y_{n+1}(q_c)\neq 0$ by our earlier remark, and since $q_c$ is
not a root of unity, we have that $\overline\y_{n+k}(q_c)\neq 0$. 

By iterating backwards one can similarly show that $\overline\y_{n-k}(q_c) \neq
0$. In fact, one can show that $\overline\y_{n+k}(q_c) + \overline\y_{n-k}(q_c) = 0$.
\end{proof}

Next we show that any such zero is a pole of $Q\left(1, 0;
\frac{q}{1+q^2}\right)$, and hence of $Q\left(1, 1;\frac{q}{1+q^2}\right)$
\begin{lemma}\label{thm:simplepole}
  Let $q = q_c \neq 0$ be a zero of $\overline\y_n$.  The function
  $Q\left(1,0;\frac{q}{1+q^2}\right)$ has a pole at the $q = q_c$.
\end{lemma}


\begin{proof}
  First, recall that as $\y_k = q^k + O(q^{k+1})$, this is a
  convergent expression (as a formal power series), and thus it makes
  sense to consider this series in this arrangement.  We contend that
  the only poles come from the poles given by zeroes of
  $\overline\y_n$, for all~$n$.  Aside from these poles and the
  square root singularity given by $q = \pm\frac{1}{\sqrt{8}}$, the
  series is convergent.

Write $\y_k\y_{k+1}$ as
\begin{equation}\label{eqnyk}
  \y_k\y_{k+1}=
  \bigl(1-q^2\bigr)\bigl(1-\y_1q\bigr)\left(\frac {1}{q^{2n}(\y_1-q) -\y_1+q}
    - \frac {1}{q^{2n+1}(q\y_1-1)-\y_1+q }
  \right).
\end{equation}
Suppose~$q$ is such that the denominator is not zero, and such
that it does not lie on the branch cut of~$\y_1$.

We know that $|q_c| \neq 1$, so we only need to consider $|q|>1$ and
$|q|<1$. The above expression for $\y_k \y_{k+1}$ is obtained by
substituting $t \mapsto \frac{q}{1+q^2}$ and so is necessarily
invariant under $q \mapsto 1/q$.  Hence we only need to consider
$|q|>1$. In this case the denominators in the right-hand side of the
expression are dominated by the $q^k$ term and so approach $0$ as $k
\to \infty$. Since the sum is alternating it converges.

Now, the denominator is annihilated in this case precisely when~$q$ is
a zero of either~$\overline\y_{k}(q)$ or $\overline\y_{k+1}(q)$. Suppose~$q_c$ is
a zero of~$\overline\y_{k}(q)$. Then, by Lemma~\ref{thm:work},
$\overline\y_{n}(q_c)\neq 0$ if $n\neq k$. Thus, the only two terms for which~$q_c$ is
singular are $\y_{k-1}\y_{k}$ and $\y_k\y_{k+1}$. If we sum these two terms
we have $(-1)^kY_{k}(Y_{k+1}-Y_{k-1})$. By our earlier comments,
$Y_{k-1}(q_c)=-Y_{k+1}(q_c)\neq 0$, and thus this singularity does not cancel.
\end{proof}

\subsection{Proof of Theorem~\ref{thm:notDfin}}
We now have all the components in place to prove the main result. 
\begin{proof}[Proof of Theorem~\ref{thm:notDfin}]
The function $Q\left(1,1; \frac{q}{1+q^2}\right)$
has a set of poles given by the zeroes of the
$\overline\y_n$, by the preceding lemma. By Lemma~\ref{thm:work}, this is an
infinite set. Thus, $Q\left(1,1, \frac{q}{1+q^2}\right)$ is not holonomic. For a multivariate series to be holonomic, any of its algebraic
specializations must be holonomic, and as  $Q\left(1,1, \frac{q}{1+q^2}\right)$ is
an algebraic specialization of both $Q(x,y)$ and $Q(1,1)$, neither
of these two functions are holonomic either. 
\end{proof}

\begin{cor}
  \label{cor asympt1}
  The number of random walks of length $n$ with step set $\{\NE, \SE,
  \NW\}$ confined to the quarter plane is asymptotic to
  \begin{equation}
    c_n \sim \alpha 3^n + O(8^{n/2})
  \end{equation}
  where $\alpha$ is a constant given by
  \begin{equation}
    \alpha =1 - 2 \sum_{n \geq 0} \frac{(-1)^n}{F_{2n} F_{2n+2}} = 0.1731788836\cdots
  \end{equation}
\end{cor}
\begin{proof}
  We need to show that the dominant singularity is the simple pole at
  $t=1/3$.  Let us rewrite equation~\eqref{eqn:kform} at $x=y=1$:
  \begin{equation}
    (1-3t) Q(1,1) = 1 - 2t Q(1,0) \qquad \mbox{ or } \qquad Q(1,1) =
    \frac{1 - 2 Q(1,0)}{1-3t}.
  \end{equation}
  Hence there are 2 sources of singularities for $Q(1,1)$ --- the
  simple pole at $t=1/3$ and the singularities of $Q(1,0)$. The
  residue at the simple pole may be computed quite directly (see the
  proof of Lemma~\ref{lem simple pole}).

  We now bound the singularities of $Q(1,0)$ away from $|t|=1/3$. The series
  $Q(1,0)$ enumerates random walks confined to the quarter plane which
  end on the line $y=0$. Consider a family of random walks with the
  same step-set but ending on the line $y=0$ but no longer confined to
  the quarter plane, rather just confined to the half-plane $y \geq 0$.
  The generating function of these walks obeys the functional equation
  \begin{equation}
    P(y)  = 1 + t\left(2y+\frac{1}{y}\right)P(y) - t\frac{1}{y}P(0).
  \end{equation}
  The generating function of these walks can be computed using the
  kernel method to give:
  \begin{equation}
    P(0) = \frac{1-\sqrt{1-8t^2}}{4t^2}.
  \end{equation}
  Hence the number of these walks grows as $O(8^{n/2})$. Consequently
  the walks enumerated by $Q(1,0)$ cannot grow faster than
  $O(8^{n/2})$, and so the radius of convergence of $Q(1,0)$ is
  bounded below by $8^{-1/2}$. The result follows.
\end{proof}

\section{The non-holonomy of  $Q_\setb(x,y;t)$}
\label{sec:setb}
This step set is not symmetric across the line $x=y$, and thus as we
lose $x \leftrightarrow y$ symmetry in the complete generating function,
we expect a more complicated scenario.  Nonetheless, we can recycle a
good number of the results from the previous case.  

\subsection{Defining two sets of iterates}
We begin by rewriting the fundamental
equation in kernel form:
\begin{equation}\label{eqn:kform2}
  K(x,y) Q(x,y) = xy - tx Q(x,0)- {ty^2}Q(0,y),
\end{equation}
with kernel $K(x,y)= xy - t(xy^2 + x^2 + y^2)$. We
will require zeros of the kernel for both fixed $x$ and fixed $y$,
which we denote respectively with $\x$ and $\y$. This gives, 
\begin{align}
  \x_{\pm 1}(y) & = \frac{y}{2t} \left(1 \mp
    \sqrt{(1-ty)^2-t^2}\right),\text{ and}\\
  \y_{\pm 1}(x) & = \frac{x}{2t(1+x)} \left( 1\mp \sqrt{1-4t^2(x+1)}\right).
 \end{align}
Using similar argument to those in Section~\ref{sec:seta}, one can show that
\begin{equation}  
\x_{\pm} ( \y_{\mp}(x) ) = x \qquad \mbox{ and } \qquad 
\y_{\pm} ( \x_{\mp}(y) ) = y
\end{equation}
and that
\begin{align}
  \frac{1}{\x_{+1}(y)} + \frac{1}{\x_{-1}(y)}
  & = \frac{1}{ty} -1  \\
  \frac{1}{\y_{+1}(x)} + \frac{1}{\y_{-1}(x)}
  & = \frac{1}{tx}.
\end{align}
We can use these to find nice expressions for compositions of
$\x$'s and $\y$'s. Define $\bX_n(y)$ and $\bY_n(x)$ by
\begin{align}
  \bX_0(y) &=y \nn
  \bX_{2k+1}(y) & = \x_1(\bX_{2k}(y))\\
  \bX_{2k+2}(y) & = \y_1(\bX_{2k+1}(y))\nonumber
\end{align}
and
\begin{align}
  \bY_0(x) &=x \nn
  \bY_{2k+1}(x) & = \y_1(\bY_{2k}(x))\\
  \bY_{2k+2}(x) & = \x_1(\bY_{2k+1}(x)).\nonumber
\end{align}
As before, we can compose these relations and determine the
recurrences
\begin{align}\label{eqn:rec2}
  \frac{1}{\bX_{n+1}(y)} + \frac{1}{\bX_{n-1}(y)}
  & = \frac{1}{t\bX_n(y)} -1  \nn
  \frac{1}{\bY_{n+1}(x)} + \frac{1}{\bY_{n-1}(x)}
  & = \frac{1}{t\bY_n(x)}.
\end{align}
Note that $\bY_k$ satisfies the same recurrence and initial conditions
as $Y_k$ in Lemma~\ref{lem kern roots seta} of Section~\ref{sec:seta}

\subsection{Expressing the series in terms of the iterates}
Next, we use the relations such as $K(\bY_{2k-1}(x), \bY_{2k}(x))=0$ and
$K(\bY_{2k}(x), \bY_{2k+1}(x))=0$ to generate an infinite set of
equations. For example, if we set $y=\bY_1(x)=\y(x)$, the functional
equation Eq.~\eqref{eqn:kform2} reduces to
\begin{equation*}
  Q(x,0) = \frac{\y_1(x)}{xt} - \left( \frac{\y_1(x)}{x} \right)^2 Q(0,\y_1(x)),
\end{equation*}
and when $x = \bX_1(y)=\x_1(y)$
\begin{equation*}
  Q(0,y) = \frac{\x_1(y)}{yt} - \left( \frac{\x_1(y)}{y} \right)^2 Q(\x_1(y),0).
\end{equation*}
Putting $y = \y_1(x)$ into the above expression for~$Q(0,y)$ gives
\begin{equation*}
  Q(0,\y_1(x)) = \frac{\x_1(\y_1(x))}{\y_1(x)t} - \left( \frac{\x_1(\y_1(x))}{\y_1(x)} \right)^2 Q(\x_1(\y_1(x)),0).
\end{equation*}
Rewriting everything in terms of $\bY_n(x)$ and substituting
one equation into the other we get
\begin{align*}
Q(x,0) &= 
\frac{\bY_0(x) \bY_1(x)}{\bY_0(x)^2 } \frac{1}{t}
- \left(\frac{\bY_1(x)}{\bY_0(x)}\right)^2
\left(
\frac{\bY_2(x)}{\bY_1(x)} \frac{1}{t} - \left( \frac{\bY_2(x)}{\bY_1(x)} \right)^2 Q(\bY_2(x),0)
 \right) \nn
& = \frac{\bY_0(x) \bY_1(x) - \bY_1(x) \bY_2(x)}{\bY_0(x)^2 } \frac{1}{t}
+ \left(\frac{\bY_2(x)}{\bY_0(x)}\right)^2 Q(\bY_2(x),0).
\end{align*}
Iterating this procedure leads to the following theorem

\begin{thm}
  \label{thm Q11 asymm}
  The generating function for the walks that end on either axis and
  the counting generating function are given by the following expressions:
  \begin{align}
    Q(x,0) &= \frac{1}{x^2 t} \sum_{n \geq 0} (-1)^n \bY_n(x) \bY_{n+1}(x), \nn
    Q(0,y) &= \frac{1}{y^2 t} \sum_{n \geq 0} (-1)^n \bX_n(y) \bX_{n+1}(y), \\
    Q(1,1) &=\frac{1}{1-3t}\left( 1 -  \sum_{n \geq 0} (-1)^n \bY_n(1) \bY_{n+1}(1) -  \sum_{n \geq 0} (-1)^n \bX_n(1) \bX_{n+1}(1)\right).\nonumber
  \end{align}
\end{thm}

\subsection{A source of singularities}
We now show that the generating function has an infinite number of
poles, and that these are given by the $\bY_n$, and the $\bX_n$.  We
consider the series under the same $q$-transformation that we used
when we considered the other step set. We set $\obY_n(q)=\bY_n(1;
\frac{q}{1+q^2})^{-1}$ and $\obX_n(q)=\bX_n(1; \frac{q}{1+q^2})^{-1}$.
\begin{lemma}
  \label{lem an eq yn}
  The function $\obY_n(q)$ is identical to $\oy_n(q)$.
\end{lemma}

\begin{proof}
  When $x=1$, the kernels for the functional equations for both step
  sets are the same. This implies that when $x=1$, the solutions to
  the respective kernels are the same; That is, 
  $\y_{+}$ is the same for both step sets and $\y_{-}$ is
  the same for both step sets.

  Consider the recurrence satisfied by the $\obY_n$. Iterating
  Eq.~\eqref{eqn:rec2} gives
\begin{equation}
  \obY_n(q) =  \left(\frac{1-q^{2n}}{1-q^2} \right)  \frac{1}{q^{n-1}} \obY_1(q)
  -  \left(\frac{1-q^{2n-2}}{1-q^2} \right)  \frac{1}{q^{n-2}},
\end{equation}
and thus, $\obY_n(q)$ here, satisfies the same recurrence as 
$\overline\y_n(q)$ given in Eq.~\eqref{eqn:bnq}, and also the same initial
condition. Thus we conclude, $\obY_n(q)=\oy_n(q)$.   
\end{proof}
\begin{cor}
  Let $q = q_c \neq 0$ be a zero of $\obY_n(q)$.  Then,
\begin{enumerate}
\item $q_c$ is not a zero for $\obY_n(q)$ whenever $n\neq k$;
\item The function  $Q\left(1,0;\frac{q}{1+q^2}\right)$ has a pole at $q =
q_c$.
\end{enumerate}
\end{cor}
\begin{proof} 
  This follows from the previous lemma and Lemmas~\ref{thm:work}
  and~\ref{thm:simplepole}.
\end{proof}

As in the case of step set $\seta$, we have that $tQ
(1,0)=\sum(-1)^n\bY_n(1;t)\bY_n(1;t)$ has an infinite number of
singularities and so in not holonomic. Thus consequently $Q(x,y)$,
is not holonomic (since one is a specialisation of the other). Thus,
already we have proven the first part of Theorem~\ref{thm:Q11}. 

In order to prove that the counting function $W(t)=Q(1,1)$ is also
not holonomic, we must show that the singularities of~$Q(1,0)$ are
singularities of $W(t)$ and so are not cancelled out by the
singularities of~$Q(0, 1)$.  This is what we do next.

\subsection{The singularities do not cancel} 
We remark that we can deduce an expression for $\bX_n\bX_{n+1}$
similar in spirit to that of Eq.~\eqref{eqnyk}. Thus, for all points
which are not zeroes of some $\obX_n(q)$, the sum $\sum_{k\geq 0}
(-1)^n\bX_k\bX_{k+1}$ converges.

Now, the singularities for ~$Q(1,0)$ and~$Q(0, 1)$ have different
sources, and thus {\em a~priori\/} we expect that they will not
be the same. We have not been able to do this. However, it is
sufficient (and easier) to show that there is an inifite subset of
singularities that do not cancel. In particular we consider the
singularities of the $\bY_n$ that lie on the imaginary axis.

For even $n$ we will show that $\bY_n$ has poles along the imaginary axis.
Now $\bX_n$ has any pole along the imaginary axis, and thus $Q(1,1)$
has an infinite number of singularities along the imaginary axis. We
prove this technical lemma now, and with that our work is done.

\begin{lemma} The functions $\obY_n(q)$ and $\obX_n(q)$ satisfy the following  recurrences:
\begin{align}\label{eq:qrec}
\obY_n&=\left(\frac{1+q^2}{q}\right)\obY_{n-1}-\obY_{n-2},\\
\obY_1&=\frac{4q}{q^2+1+\sqrt{q^4-6q^2+1}} & \obY_0=1\nn
\obX_n&=\left(\frac{1+q^2}{q}\right)\obX_{n-1}-\obX_{n-2}-1,\nn
\obX_1&=\frac{2q}{q^2-q+1+\sqrt{q^4-2q^3-q^2-2q+1}}& \obX_0=1.\nonumber
\end{align}
Thus, it follows
\begin{itemize}
\item[(1)] For $n \geq 1$, $\obY_{2n}$ has at least two zeros along the imaginary axis;
\item[(2)] For any $n$, and any $q=i r$ along the imaginary axis the
  value of $\obX_n$ is non-zero, except at
  $q=0$. 
\end{itemize}
\end{lemma}

\begin{proof} Equation~\eqref{eq:qrec} is a direct consequence of Equation~\eqref{eqn:rec2}.
\noindent{(1)} From the proof of Lemma~\ref{thm:work} in
Section~\ref{sec:seta} we see that the zeros of $\overline{Y}_n$ are
given by the zeros of the following equation
\begin{equation}
    (1-q^{2n}) \sqrt{1-6q^2+q^4} + (1-q^2)(1+q^{2n}) = 0.
\end{equation}
By Lemma~\ref{lem an eq yn} we know that the zeros of
$\obY_n$ satisfy the same equation.

Now substitute $q = ir$ and $n = 2m$ into this equation to obtain
\begin{equation}
  (1-r^{4m}) \sqrt{1+6r^2+r^4} + (1+r^2)(1+r^{4m}) = 0.
\end{equation}
Denote the left-hand side of this equation by $f(r)$, which is a
continuous function for real $r$.

At $r = \pm 1$ we find that $f(\pm 1) = 4$. At $r = \pm 2$ we find
that
\begin{equation}
  f(\pm 2)  = (\sqrt{41}+5) + (5-\sqrt{41}) 2^{4m} 
\end{equation}
which is negative for $m \geq 1$.  Hence $f(r)$ has at least two
zeros, one in the interval $[-2,-1]$ and the other in the interval
$[1,2]$. So $\overline{A}_{2n}$ has zeros on the imaginary axis.

\noindent{(2)} First, we solve the recurrence and give an explicit
form for~$\overline{\bX}_n$ in terms of $\overline{\bX}_1$, denoted
here by $\beta$:
\begin{equation}\label{eqn:explicit}
\overline{\bX}_n=\frac{q^{2n} (1-2q+\beta q^2-\beta q) +q^n(q^2+q)
+(q^3-2q^2+q\beta-\beta q^2)}%
{ \left( -1+q \right) ^{2} \left( q+1 \right) q^n}.
\end{equation}

We can now repeat the method used in the proof of Lemma~\ref{thm:work}
to show that the zeros of $\obX_n$ satisfy the following equation:
\begin{multline}
  0 = q^2 (q^2-3q+1) \left(q^{4n}+1\right) 
  - q (q+1)^2(q^2-3q+1)\left(q^{3n}+q^n \right)
  \\ - (q^6-4q^5+3q^4+6q^3+3q^2-4q+1) q^{2n}.
\end{multline}
Note that this equation is symmetric in $q \leftrightarrow 1/q$ (as we
would expect, since we introduced the $q$ variable by substituting $t
\mapsto 1/(q+1/q)$). 

We wish to show that the solutions of this equation do not lie on the
imaginary axis and so the singularities of $\bY_n$ which lie on the
imaginary axis cannot be cancelled by singularities of $\bX_n$. Hence
we substitute $q = i r$ with $r \in \mathbb{R}$:
\begin{multline}
  \label{eqn nu arr}
  0 = r^2 (r^2+3ir-1) \left(r^{4n}+1\right) 
  - i r (r-i)^2(r^2+3ir-1)\left((ir)^{3n}+(ir)^n \right)
  \\ + (r^6+4ir^5-3r^4+6ir^3+3r^2+4ir-1) r^{2n} (-1)^n.  
\end{multline}
Let $\nu(r)$ denote the right-hand side of this expression. Note that
$\nu(0)=0$ for all $n$ and also $\nu(\pm1) = \pm8I$, so that $r=0$ is
always a zero and $r = \pm1$ are never zeros.

We now analyse the real zeros of this expression for different values
of $n \mod 4$.  In each case we show that there are no real values for
$r$ which annihilate~$\nu$ (except $r=0$). We do this by considering
the real and imaginary parts of $\nu$ separately.

\subsection*{Case 1:  $n\equiv 2\pmod 4$}
In the case that $n\equiv 2\pmod 4$, the imaginary part of $\nu$ is 
\begin{equation}
  \Im(\nu)= 
3r^{4n+3}+r( r^4 + 4r^{2}+1) {r}^{3n}
+2r( 2r^4+3r^{2}+2) {r}^{2n}
+r( r^4 + 4r^2+1)r^n+3r^3.
\end{equation}
This expression is strictly positive for
$r>0$. Since it contains all odd powers of $r$, it is strictly
negative for $r<0$. It can only be zero when $r=0$.

\subsection*{Case 2: $n\equiv 0\pmod 4$}
In this case, $\nu$ has real part equal to
\begin{equation}
  \Re(\nu)=r^2(r^2-1)\left(r^{4n}+r^{3n}+(r^2-2+r^{-2}) r^{2n}+r^n+1\right).
\end{equation}
Suppose $r \neq 0, \pm 1$.  Then
\begin{equation}
  r^2\left(r^{4n}+r^{3n}+r^{n}+r + (r^2-1)^2  r^{n}\right)=0
\end{equation}
implies that $r^{2n}+r^n +r^2+r^{-2} + r^{-n}+r^{-2n}=2$. This
equation is symmetric in $r \leftrightarrow 1/r$ and $r
\leftrightarrow -r$ (since $n$ is even). If $r \geq 1$, the left-hand
side is strictly larger than 2, and hence has no solution. By
symmetries, there can be no real solution.

So the only points that annihilate the real part of $\nu$ are $r=0$
and $r = \pm1$. As noted above $r=0$ is a zero of $\nu$, but $r = \pm
1$ are not.

\subsection*{Case 3: $n$ odd} To show that~$\nu$ has no real zeros
(aside from $r=0$) when $n$ is odd, we consider a linear combination
of the real and imaginary parts of~$\nu$ that eliminates the leading
powers of $r$ which in turn simplifies the analysis. 

In particular consider $(r-1/r) \Im(\nu) - 3\Re(\nu)$; this eliminates
the $(r^2+3ir-1)(r^{4n}+1)$ term from equation~\eqref{eqn nu
  arr}. If there exists a real zero of $\nu$ then it must also
annihilate this linear combination. We will show that no such value of
$r$ exists, and so there are no real zeros of $\nu$ (excepting $r=0$).

The linear combination simplifies to
\begin{multline}
(r-1/r) \Im(\nu)-3\Re(\nu)\\
 = r^{2k+2}(2+14r^2+2r^4)(r^{4k+2}-1)(-1)^k - (r^2-1) (r^4+12r^2+1)r^{4k+2}.
\end{multline}
Now divide both sides by $r^2-1$ to obtain (we know that $r \neq \pm
1$ are not zeros of $\nu$):
\begin{multline}
  \frac{1}{r}\Im(\nu)-\frac{3}{r^2-1}\Re(\nu)\\
   =
  r^{2k+2}(2+14r^2+2r^4)\left(\frac{r^{4k+2}-1}{r^2-1}\right)(-1)^k 
  - (r^4+12r^2+1)r^{4k+2}.
\end{multline}
When $k$ is odd this expression is a polynomial in $r^2$ with
negative coefficients, while when $k$ is even this is a polynomial in
$r^2$ with positive coefficients. Hence this is non-zero except at $r=0$.

Thus, as no $\overline{\bX}_n$ has a zero on the imaginary axis, (2)
holds.
\end{proof}

It seems clear that  $Q(0,1;t)$ is also non-holonomic. Indeed, the
walks are very similar to the symmetric case. However, we were unable
to construct a direct argument to show that it too is non-holonomic.

\section{Conclusions}
\subsection{How robust is this method at detecting non-holonomy?}
This is a natural question, and it speaks to the ability of this
approach to be applied to other problems. In all the cases known to
give a holonomic generating function, the group of kernel iterates is
finite (of order 8 or less). If one attempts to repeat the proofs of
Theorems~\ref{thm Q11 symm} and~\ref{thm Q11 asymm} then one obtains a
finite set equations. It is not clear that solution of these equations
has singularities appearing in the same way as occurred here.

\subsection{Combinatorial intuitions of non-holonomy}
More generally, we are interested in developing an intuition about
when precisely a combinatorial object will have a non-holonomic
generating function. If we decompose these walks according to their
$\NE$ or $\NN$ steps, a source of singularities becomes
apparent. Unfortunately, since we are unable to prove that the
singularities do not cancel, the details we present in this section do
not constitute a proof of non-holonomy, however they are nonetheless
instructive for understanding a potential source of the singularities.

\begin{figure}
\center
\includegraphics[scale=0.8]{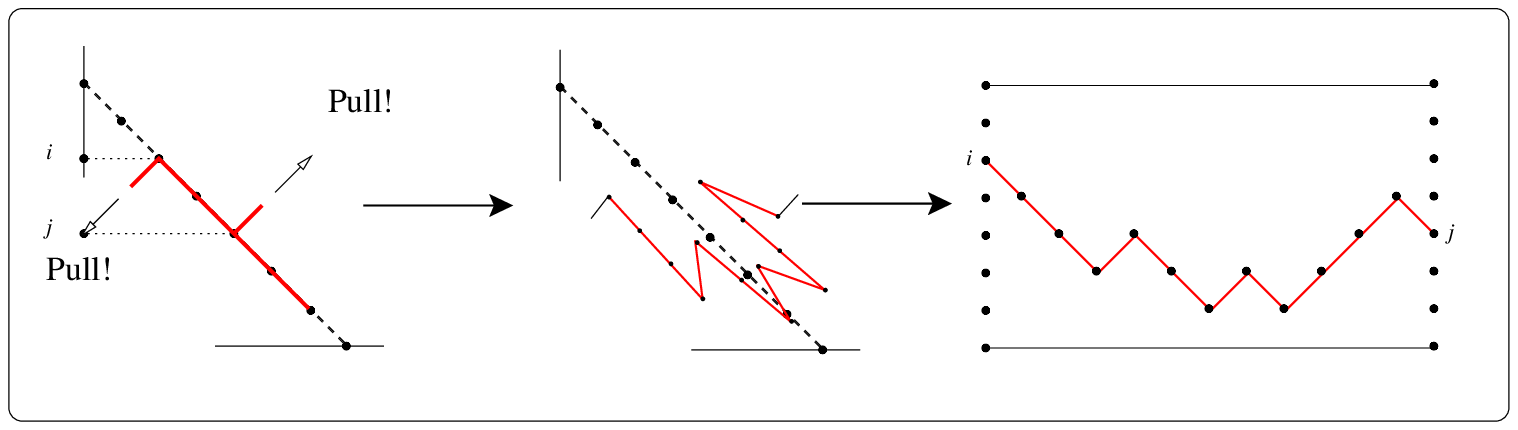}
\caption{Stretching the walk to find a directed path in a strip}
\label{fig:stretch}
\end{figure}

Any walk from either of these two families can be decomposed into a
sequence of directed paths in a strip of unrestricted length where after
each $\NE$/$\NN$ step the size of the height bound
increases. Figure~\ref{fig:stretch} for an example of such a
decomposition for step set~$\seta$.
Thus, a potential strategy groups lattice walks that end on the line
$x+y=k$,  and decomposes any such walk into a triple: A walk that ends on the line
$x+y=k-2$, a $\NE$ step, and a directed path in a strip of height $k$.

To describe the generating function, we use Example~11 of~\cite{FlSe93} who
present a generating function for paths of length $n$ in strip of
height $k$ of walks that begin at a given height, and end at a given height.

Define $D_k(y;t)$ as
the generating function for the subset of walks ending on $x+y=k$
where $y$ marks the final height of the walk. We can easily translate
the above decomposition into a functional equation for the generating
function. If we then make the substitution $t\mapsto\frac{q}{1+q^2}$,
the expression simplifies remarkably into the following recurrence for $D_k(y)=D_k(y, \frac{q}{1+q^2})$:
%
\begin{equation}
D_k(y)=\frac{q^3 D_{k-2}(q)(y^{k+2}+1)- qy^2D_{k-2}(y)(q^{k+2}+1)}
{(q^{k+2}+1)(yq-1)(y-q) }.
\end{equation}
In fact, for our purposes it suffices to consider:
\begin{equation}
D_k(1)=\frac{q(q^{k+2}+1)D_{k-2}(1)-2q^3 D_{k-2}(q)}{(q^{k+2}+1)(q-1)^2}.
\end{equation}
From this formula, and from computations for various values of $k$,
$D_k(1)$ is a rational function in $q$, and it seems clear that the
set (taken over all $k$) of poles of $D_k(1)$ is dense in the unit
circle. Were this so, we would apply the following theorem to the
generating function $Q_\seta(s, s; \frac{q}{1+q^2})=\sum
D_k(1;\frac{q}{1+q^2})s^k$, and thus conclude the non-holonomy of
$Q_\seta(x,y;t)$.
\begin{thm}
  Let $f(x;t)=\sum_n c_n(x) t^n$ be a holonomic power series in $\bC(x)
  [[t]]$. with rational coefficients in $x$. For $n\geq 0$ let $S_n$
  be the set of poles of $c_n(y)$, and let $S=\bigcup S_n$. Then $S$
  has only a finite number of accumulation points.
\end{thm} 
Again, the principle difficulty is showing that the singularities do
not cancel; that solutions to $q^{k+2}+1$ are indeed poles of $D_k$. 

This approach was pioneered by Guttman and Enting~\cite{GuEn96}, and
has been fruitful for several different models~\cite{Guttman00,Rechnitzer06}. Unfortunately their arguments do not appear to work
here. 
\subsection{Related walks}
We expect the following walks to have non-holonomic generating
functions because the groups of their kernel iterates are infinite. It
is even likely this can be proved in the same manner as Theorems~\ref{thm:notDfin} and~\ref{thm:Q11}.
\begin{center}
\includegraphics[width=5cm]{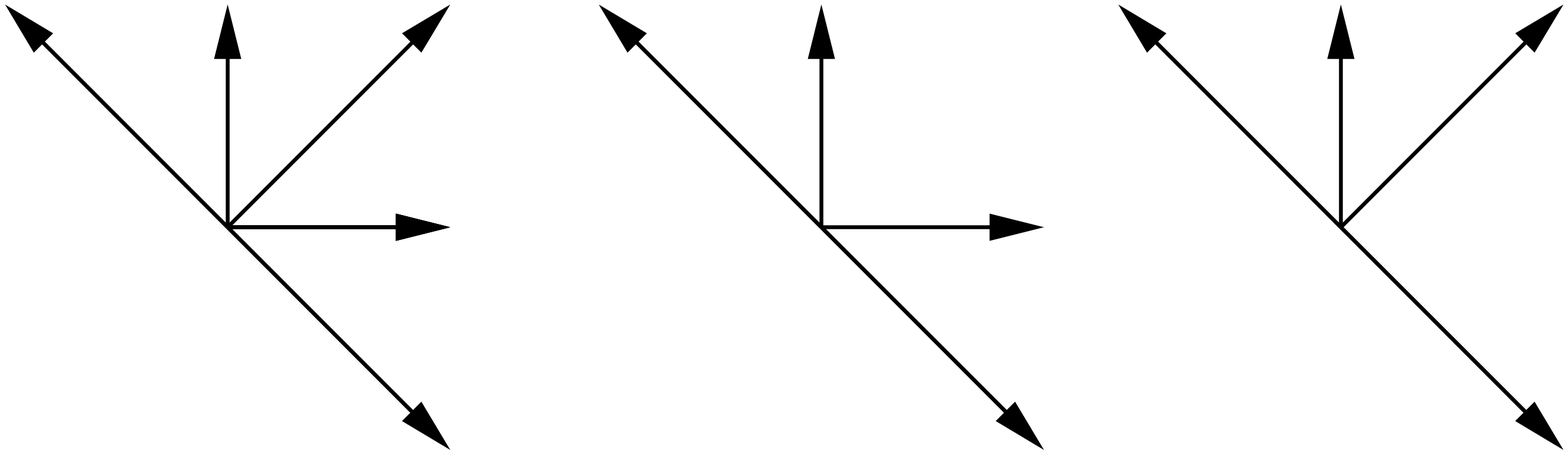}
\end{center}
\subsection{Acknowledgments}
We wish to thank Mireille Bousquet-M\'elou and Bruno Salvy for
extremely useful discussions. Indeed, some of this work was completed
during a postdoc of the first author with Bousquet-M\'elou at LaBRI, Universit\'e Bordeaux
I. Cedric Chauve offered many useful comments to improve
the readability. This work was funded in part by an NSERC Discovery
grant, and an NSERC Postdoctoral Fellowship.

\bibliographystyle{plain}

\begin{thebibliography}{10}

\bibitem{BaFl02}
Cyril Banderier and Philippe Flajolet.
\newblock Basic analytic combinatorics of directed lattice paths.
\newblock {\em Theoret. Comput. Sci.}, 281(1-2):37--80, 2002.
\newblock Selected papers in honour of Maurice Nivat.

\bibitem{MBM05}
Mireille Bousquet-Melou.
\newblock {Walks in the quarter plane: Kreweras' algebraic model}.
\newblock {\em Annals of Applied Probability}, 15(2):1451--1491, 2005.

\bibitem{MBM06}
Mireille Bousquet-M\'elou.
\newblock Rational and algebraic series in combinatorial enumeration.
\newblock In {\em International Congress of Mathematicians 2006}, 2006.

\bibitem{BoPe03}
Mireille Bousquet-M{\'e}lou and Marko Petkov{\v{s}}ek.
\newblock Walks confined in a quadrant are not always {D}-finite.
\newblock {\em Theoret. Comput. Sci.}, 307(2):257--276, 2003.
\newblock Random generation of combinatorial objects and bijective
  combinatorics.

\bibitem{FlSe93}
Philippe Flajolet and Robert Sedgewick.
\newblock Analytic combinatorics: functional equations, rational and algebraic
  functions.
\newblock {\em INRIA Research Report}, 4103, January 2001.

\bibitem{Guttman00}
A.J. Guttmann.
\newblock Indicators of solvability for lattice models.
\newblock {\em Discr. Math.}, 217:167--189, 2000.

\bibitem{GuEn96}
A.J. Guttmann and I.G. Enting.
\newblock On the solvability of some statistical mechanical systems.
\newblock {\em Phys. Rev. Letters}, 76:344--347, 1996.

\bibitem{Lipshitz89}
Leonard Lipshitz.
\newblock {$D$}-finite power series.
\newblock {\em J. Algebra}, 122(2):353--373, 1989.

\bibitem{MishnaXX}
Marni Mishna.
\newblock Classifying lattice walks in the quarter plane.
\newblock {\em preprint: ArXiv:math.CO/0611651}, 2006.

\bibitem{Prodinger03}
Helmut Prodinger.
\newblock The kernel method: a collection of examples.
\newblock {\em S\'em. Lothar. Combin.}, 50:Art. B50f, 19 pp. (electronic),
  2003/04.

\bibitem{JaPeRe06}
T.~Prellberg, A.~Rechnitzer, and E.~J.~Janse van Rensburg.
\newblock Partially directed paths in a wedge.
\newblock {\em preprint arxiv.org:math.CO/0609834}, 

\bibitem{Rechnitzer06}
A.~Rechnitzer.
\newblock Haruspicy 2: The anisotropic generating function of self-avoiding
  polygons is not d-finite.
\newblock {\em J. Combinatorial Theory - Series A}, 113:520--546, 2006.


\bibitem{Stanley80}
R.~P. Stanley.
\newblock Differentiably finite power series.
\newblock {\em European J. Combin.}, 1(2):175--188, 1980.
\end{thebibliography}

\end{document}